\documentclass[12pt]{amsart}

\usepackage{ulem}

\usepackage{amsmath,amssymb,amsthm,enumitem}

\usepackage{multirow}

\usepackage{tikz-cd}

\usepackage{amsfonts}

\usepackage{tikz}
\usetikzlibrary{matrix,arrows,decorations.pathmorphing}

\usepackage{enumitem}

\makeatletter
\@namedef{subjclassname@2010}{%
  \textup{2010} Mathematics Subject Classification}
\makeatother

\DeclareMathOperator{\Supp}{Supp}

\theoremstyle{definition}

\newtheorem*{Lemule}{Lemma}

\def\cqfd{
{\hfill
\kern 6pt\penalty 500
\raise -1pt\hbox{\vrule\vbox to 5pt{\hrule width 4pt
\vfill\hrule}\vrule}}
\break}

\font\tengoth=eufm10
\font\sevengoth=eufm7
\font\fivegoth=eufm5
\newfam\gothfam
\textfont\gothfam=\tengoth\scriptfont\gothfam=\sevengoth\scriptscriptfont\gothfam=\fivegoth

\title[Hypersurfaces in weighted projective spaces]{Corrigendum to the paper ``Maximum number of rational points on hypersurfaces in weighted projective spaces over finite fields'''\\
Journal of Algebra and Its Applications, Vol. 24, No. 13n14, 2541015
(2025)
World Scientific Publishing Company
DOI: 10.1142/S0219498825410154}

\author[Yves Aubry]{Yves Aubry}
\address[Aubry]{Institut de Math\'ematiques de Toulon - IMATH, Universit\'e de Toulon, France}
\address[Aubry]{Institut de Math\'ematiques de Marseille - I2M, Aix Marseille Univ, UMR 7373 CNRS, France}
\email{yves.aubry@univ-tln.fr, yves.aubry@univ-amu.fr}

\author[Marc Perret]{Marc Perret}
\address[Perret]{Institut de Math\'ematiques de Toulouse - UMR 5219, CNRS, UT2J, F-31058 Toulouse, France}
\email{perret@math.univ-toulouse.fr}

\date{\today}

\begin{document} 

\baselineskip=17pt



\maketitle

\begin{abstract}
The statement of item (ii) of  Proposition 3.2 of the article referenced in the title
is not correct. We provide a corrected version and show that, under the assumption
that $\gcd(a_i, a_j, q-1)=1$ for any pair $i\neq j$ in $\{0, \cdots, n\}$ (with the notations of the paper),
our initial statement becomes valid, as does the remainder of the paper.

\end{abstract}

\bigskip

As pointed to us by Jade Nardi and 
Rodrigo San-Jos\'e,
Proposition 3.2 (ii) of our paper cited in the title is not correct.
More precisely, let $P$ be a point lying in the set 
$${\mathcal T}_i=\left\{P=[y_0:\cdots:y_n] \in {\mathbb P}(a_0, \cdots, a_n)({\mathbb F}_q)  \vert y_i=1 ;  P\neq [0:\cdots 0:1:0\cdots :0]\right\}.$$ 
Then, the number $\sharp \pi_i^{-1}(P)({\mathbb F}_q)$ of ${\mathbb F}_q$-rational preimages of $P$ by 
$$\begin{matrix} \pi_i & : & {\mathbb P}(a_0, \cdots, a_{i-1}, 1, a_{i+1}, \cdots, a_n)  &\rightarrow &{\mathbb P}(a_0, \cdots, a_n)\\&&[x_0: \cdots:x_n] & \mapsto&[x_0: \cdots: x_{i-1}:x_i^{a_i}:x_{i+1}:\cdots:x_n]\end{matrix}$$
is not equal to
$$\gcd(a_i,q-1)$$
as claimed in our Proposition~3.2. Indeed, let us for instance consider, for fixed $a_0, a_1 \geq 1$, the case
$$\begin{matrix}\pi_2 : & {\mathbb P}(a_0,a_1, 1,4)& \rightarrow & {\mathbb P}(a_0,a_1, 2,4)\\&[x_0:x_1:x_2:x_3]&\mapsto & [x_0:x_1:x_2^2:x_3]\end{matrix}$$
for $q=5$, and the point $P =[0,0,1,2] \in  {\mathcal T}_2 \subset {\mathbb P}(a_0,a_1, 2,4)({\mathbb F}_5)$. We thus have $a_2=2$ and $q-1=4$, so that our Proposition~3.2 predicts $\gcd(2,4)=2$ rational inverse preimages $Q$ of $P$ in  ${\mathbb P}( a_0,a_1, 1,4)({\mathbb F}_5)$, while a close study shows that it has only one, namely $Q=[0:0:1:2]=[0:0:-1:2]$. 

\medskip

Instead, the following Lemma is true.

\begin{Lemule} [Corrected form of Item (ii) of Proposition 3.2]
Let $i\in\{0,\ldots,n\}$
and $P=[y_0\cdots:y_n] \in {\mathcal T}_i \subset {\mathbb P}(a_0, \cdots, a_n)({\mathbb F}_q)$ with $y_i=1$ and $y_j \in {\mathbb F}_q$ for $0 \leq j \leq n$. 

Let $\delta_{P}=\gcd(a_j \vert   j\in \Supp(P))$ and  $\delta_{i,P}=\gcd(a_j \vert  j\in \Supp(P) \hbox{~and~} j\neq i)$
where  $\Supp(P)$ denotes the support of $P$.

\begin{enumerate}
\item \label{formule} Then, we have
$$\sharp \pi_i^{-1}(P)({\mathbb F}_q)=\frac{\gcd(a_i,(q-1)\times \delta_{i, P})}{\delta_{P}}.$$

\item \label{Hyp} Assuming that 
$\gcd(a_i, a_j, q-1)=1$
for any $j\in\{1,\ldots,n\}\setminus \{i\}$, 
this reduces to 
$$\sharp \pi_i^{-1}(P)({\mathbb F}_q) = \gcd(a_i,q-1).$$
\end{enumerate}
\end{Lemule}

Note that the assumption in 
Item~(\ref{Hyp}) 
is trivially satisfied, for instance
\begin{itemize}
\item either if the weights $a_i, a_j$ are coprime for any $j\neq i$, 
\item or if $a_i$ is coprime to $q-1$.
\end{itemize}

Under one of the above extra conditions, the result 
stated in Proposition~3.2~(ii) is thus correct, as well as the whole paper.

\bigskip

{\it Proof of the Lemma.} Let us assume for convenience that $i=0$, and that
$$P=[1 : y_1: \cdots : y_m :0:\cdots:0]\in {\mathcal T}_0\subset {\mathbb P}(a_0, a_1, \cdots, a_n)({\mathbb F}_q),$$
with $y_j \in {\mathbb F}_q$ and $y_j \neq 0$ for $1\leq j\leq m$. Note that 
we have necessarily $m\geq 1$
since $[1:0:\cdots\:0]\notin {\mathcal T}_0$.

We begin by describing the whole set $\pi_i^{-1}(P)(\overline{\mathbb F}_q)$. Let $Q = [x_0:\cdots:x_n] \in  {\mathbb P}(1, a_1, \cdots, a_n)(\overline{\mathbb F}_q)$
and 
let us denote by 
$\mu_r(\overline{\mathbb F}_q)$ the set of $r$-th roots of unity in $\overline{\mathbb F}_q$.

We have $Q \in \pi_i^{-1}(P)$ if and only if
$[x_0^{a_0} : x_1 : \cdots : x_n]=[1 : y_1: \cdots : y_m :0:\cdots:0]$  in ${\mathbb P}(a_0, a_1, \cdots, a_n)$
which means  that 
there exists $\lambda \in \overline{\mathbb F}_q^*$ such that
$$\left\{ \begin{matrix} x_0^{a_0}&=&\lambda^{a_0}\times 1&\\x_j&=&\lambda^{a_j} \times y_j & \quad (1\leq j\leq m)\\ x_j&=&\lambda^{a_j} \times 0 & \quad (m<j).\end{matrix} \right.$$

This is equivalent to saying that 
there exists $\lambda \in \overline{\mathbb F}_q^*$ and  $\zeta \in \mu_{a_0}(\overline{\mathbb F}_q)$ such that
$$\left\{ \begin{matrix} x_0&=&\lambda \times \zeta&\\x_j&=&\lambda^{a_j} \times y_j & \quad (1\leq j\leq m)\\ x_j&=& 0 & \quad (m<j)\end{matrix} \right.$$
i.e. 
such that
$Q=[x_0:\cdots:x_n] = [\lambda\times \zeta : \lambda^{a_1}\times y_1: \cdots :  \lambda^{a_m}\times y_m: 0\cdots 0].$

Thus we have proved that
$$\pi_i^{-1}(P)(\overline{\mathbb F}_q) = \left\{Q_{\zeta}=[ \zeta : y_1: \cdots :  y_m: 0\cdots 0] ; \quad \zeta \in \mu_{a_0}(\overline{\mathbb F}_q)\right\}.$$

Next, we determine the ${\mathbb F}_q$-rational points inside the above set 
$\pi_i^{-1}(P)(\overline{\mathbb F}_q)$.
Let $\zeta \in \mu_{a_0}(\overline{\mathbb F}_q)$. 
From $\zeta \neq 0$ and $y_i^q=y_i\neq 0$ for $1\leq i\leq m$ (since $y_i \in {\mathbb F}_q^*$), we have that 
$Q_{\zeta}=[ \zeta : y_1: \cdots :  y_m: 0\cdots 0] \in \pi_i^{-1}(P)({\mathbb F}_q)$
if and only if $[ \zeta^q : y_1^q: \cdots :  y_m^q: 0^q\cdots 0^q]=[ \zeta : y_1: \cdots :  y_m: 0\cdots 0]$
which is equivalent to saying that
there exists $\lambda \in \overline{\mathbb F}_q^*$ such that 
$$\left\{ \begin{matrix} \zeta^q&=&\lambda \times \zeta&\\y_j^q&=&\lambda^{a_j} \times y_j & \quad (1\leq j\leq m)\\ \end{matrix} \right.$$
i.e. such that 
$$\left\{ \begin{matrix} \lambda&=\zeta^{q-1}&\\\lambda^{a_j} &=1& \quad (1\leq j\leq m).\end{matrix} \right.$$
This means that $\zeta^{(q-1)a_j}=1$ for all  $1\leq j\leq m$, in other words that
$$ \zeta \in \cap_{1\leq j\leq m}\mu_{(q-1)a_j}(\overline{\mathbb F}_q)=\mu_{\gcd((q-1)a_1,\cdots, (q-1)a_m)}(\overline{\mathbb F}_q).$$
It follows that
$\pi_i^{-1}(P)({\mathbb F}_q)$ is the set of points 
$Q_{\zeta}=[ \zeta : y_1: \cdots :  y_m: 0\cdots 0]$
such that   $\zeta \in \mu_{a_0}(\overline{\mathbb F}_q)\cap\mu_{(q-1)\times \delta_{0,P}}(\overline{\mathbb F}_q)=\mu_{\gcd(a_0,(q-1)\times\delta_{0,P})}(\overline{\mathbb F}_q)$
where $\delta_{0,P}=\gcd(a_1, \cdots, a_m)$.

We now need to determine the number of distinct elements in this set.
Let $\zeta_1, \zeta_2 \in \mu_{\gcd(a_0,(q-1)\times\delta_{0,P})}(\overline{\mathbb F}_q)$. 
We have $Q_{\zeta_1}=Q_{\zeta_2}$ if and only if 
$[\zeta_1 : y_1: \cdots :  y_m: 0\cdots 0] = [\zeta_2 : y_1: \cdots :  y_m: 0\cdots 0]$  in ${\mathbb P}(1,a_1, \cdots, a_n)$.
This is equivalent to the existence of 
 $\lambda \in \overline{\mathbb F}_q^*$ such that
$$\left\{ \begin{matrix} \zeta_1&=&\lambda \times \zeta_2&\\y_j&=&\lambda^{a_j} \times y_j & \quad (1\leq j\leq m)\\ \end{matrix} \right.$$
i.e. such that
$$\left\{ \begin{matrix} \lambda&=&\zeta_1/\zeta_2 &\\\lambda^{a_j} &=&1 & \quad (1\leq j\leq m).\\ \end{matrix} \right.$$
It writes 
$(\zeta_1/\zeta_2)^{a_j} =1$ for all $1\leq j\leq m$, hence we have proved that 
$$Q_{\zeta_1}=Q_{\zeta_2} \iff \zeta_1/\zeta_2 \in \mu_{\gcd(a_0,(q-1)\times \delta_{0,P})}(\overline{\mathbb F}_q)\cap\mu_{\delta_{0,P}}(\overline{\mathbb F}_q)=\mu_{gcd(a_0,\delta_{0,P})}(\overline{\mathbb F}_q).$$

We deduce that 
$$\sharp \pi_i^{-1}(P)({\mathbb F}_q)=\frac{\gcd(a_0,(q-1)\times \delta_{0, P})}{\gcd(a_0, \delta_{0,P})},$$
which proves
 Item~(\ref{formule}).

\bigskip

 In order to prove 
Item~(\ref{Hyp}), 
  let $\ell$ be any prime number and let us set
$$\alpha:=v_{\ell}(a_0), \quad \kappa:=v_{\ell}(q-1) \hbox{~and~} \delta:=v_{\ell}(\delta_{0,P})$$
where $v_{\ell}$ stands for the $\ell$-adic valuation.
We have
\begin{align}\label{vell-1}
v_{\ell}\left(\frac{\gcd(a_0,(q-1)\times \delta_{0, P})}{\gcd(a_0, \delta_{0,P})}\right)
&=\min(\alpha, \kappa+\delta)-\min(\alpha, \delta)
\end{align}
while
\begin{align}\label{vell-2}
v_{\ell}(\gcd(a_0, q-1))
&=\min(\alpha, \kappa).
\end{align}

Under the extra assumption that $\gcd(a_i, a_j, q-1)=1$, at least one of the three valuations $\alpha$, $\kappa$ or $\delta$ do vanish. In each case, it is a trivial matter to observe that the right hand side in Equations~(\ref{vell-1}) and~(\ref{vell-2}) are equal, so as their left hand side which proves
 Item~(\ref{Hyp}).
 
 \bigskip
 
 We end this note by drawing attention to the  preprint ``Maximum number of zeroes of polynomials on weighted projective spaces over a finite field'', arXiv:2507.22597v1 [math.AG] 30 Jul 2025, 
 in which the authors Jade Nardi and Rodrigo San-Jos\'e 
 present a proof of the  conjecture in the general case.
 
 \bigskip

\end{document}